\input amstex
\documentstyle{amsppt}
\input boxedeps.tex
\SetepsfEPSFSpecial
\HideDisplacementBoxes
\hsize=5.5in
\vsize=7.4in
\magnification=1200
\def\e{\epsilon}

\def\i{^{-1}}
\def\l{\lambda}

\def\a{\alpha}

\def\G{\Gamma}

\def\<{\langle}
\def\>{\rangle}
\parindent=18pt

\document
\baselineskip=18pt
\NoBlackBoxes
\def\leaderfill{\leaders\hbox to 1em{\hss.\hss}\hfill}

{\centerline{\bf{GROUPS GENERATED BY POWERS OF DEHN TWISTS}}}
{\centerline{\bf{}}}

\

{\centerline{{Hessam Hamidi-Tehrani}
\footnote{Department of Mathematics, University of California at
Santa Barbara, CA 93106
E-mail: hessam\@math.ucsb.edu.
http://www.math.ucsb.edu/$\sim$hessam.}
\footnote{Partially supported by an NSF grant.}}}
\

{\centerline{{April 15, 1998}}}

\

\topmatter
\abstract
We classify groups generated by powers of 2 Dehn twists which are 1) free or
2) have no ``unexpected'' reducible elements. We give some sufficient conditions in the case of groups generated by powers of $h \ge 3$ twists.

\


\


\endabstract
\endtopmatter

\head \S 0. Introduction \endhead

Let $S$ be a punctured or non-punctured oriented surface. For (the isotopy 
class of)\footnote{ We will usually drop this in the rest of this paper for
 brevity.} a simple  closed curve $c$
on $S$ let $D_c$ denote the right-handed Dehn twist about $c$. Let 
$(c_1,c_2)$ denote the geometric intersection number of simple closed curves
 $c_1,c_2.$ Let $M(S)$ be the mapping class group of $S$. Let $\Bbb F_h$ be the free group on $h$ generators.
For a set of $h$ simple closed curves $A=\{a_1,\cdots a_h\}$ and positive 
integers 
$n_1,\cdots n_h$ we study the group $G= \< D^{n_1}_{a_1},\cdots D^{n_h}_{a_h}\>$, and ask the question whether $G \cong \Bbb F_h$.
It is well known that $G \cong \Bbb F_h$ if $(a_i,a_j) \ne 0$ for $i \ne j$ and $n_i$  are large, for all $i,j$ (See for example [I]).

In the case $h=2$ we will give a complete answer, i.e., $G \ncong \Bbb F_2$
if and only if 
$$((a_1,a_2), \{n_1,n_2 \})=(0,\{*,*\}), (1,\{1,1\}),
 (1,\{1,2\}),{\text { or }} (1, \{1,3\})$$
 (See Theorem 2.4).

It should be noticed that in the case of two curves
 $a,b$ filling up a closed surface this was done by Thurston as a method
 to construct pseudo-Anosov elements;
i.e., he showed that $\<D_a,D_b\>$ is free and consists of pseudo-Anosov 
elements besides powers of conjugates of the generators
[FLP]. Our methods are completely different and elementary, and are only based on how the geometric intersection pairing behaves under Dehn twists.

In the case when $h \ge 3$, we give some sufficient conditions for $G \cong \Bbb F_h$. To motivate our condition, look at the case
 $\G=\<D_{a_1}, D_{a_2}, D_{a_3} \>$ , and assume
$a_3=D_{a_1}(a_2)$. Now $D_{a_3}=D_{a_1}D_{a_2}D_{a_1}\i$, so $G \ncong \Bbb F_3$. But notice that  $(a_1,a_3)=(a_1,a_2)$ and $(a_2,a_3)=(a_1,a_2)^2$.
This shows that the set $I=\{(a_i,a_j) \ | \ i\ne j \}$ is ``spread around''.
It turns out this is in a sense the only obstruction for $\G \cong \Bbb F_h$:

\proclaim{Theorem 0.1} Suppose $\G= \< D_{a_1},\cdots D_{a_h}\>$,
and let $m=\min I$ and $M=\max I$, where $I=\{(a_i,a_j) \ | \ i\ne j \}$.
Then $\G \cong \Bbb F_h$ if $M \le m^2/6$.
\endproclaim

We will prove a more general version of this (see Theorem 3.2).

The second question we ask about the group $G$ is that to what extent the elements of $G$ are pseudo-Anosov? a mapping class $f$ is called pseudo-Anosov if
$f^n(c) \ne c$ for all non-trivial simple closed curves and $n>0$.
Let $A=\{a_1,a_2, \cdots, a_h\}$  be a set of simple closed curves on $S$.
The surface filled by $A$, denoted $S_A$ is a regular neighborhood $N$ of 
$\cup A$ together with the components of $S-N$ which are discs with 0 or
 1 punctures. We say that $A$ fills up $S$ if $S_A=S$.
Let $f=D_{c_1}^* \cdots D_{c_k}^*$ be a cyclically-reduced word in the generators of $G$. Let
$supp(f)=S_{\{c_1,\cdots, c_k\}}$. Then we say $f$ is relatively pseudo-Anosov if $f$ is either the identity or
pseudo-Anosov in $M(supp(f))$.  If $g=h f h\i$ with $f$ cyclically-reduced, define $supp(g)=h(supp(f))$. The group $G$ is relatively  pseudo-Anosov if every word in generators of $G$ is relatively pseudo-Anosov.
Notice that a Dehn twist $D_{a_i}$ is relatively pseudo-Anosov since 
its support is a twice-punctured sphere.
Intuitively, a group $G$ is relatively 
pseudo-Anosov if no element of $G$ has ``unexpected
reducibility''.  It is well known that for $n_i$ large, 
$G$ is relatively pseudo-Anosov for any choice of curves $a_i$.

In the case $h=2$ again we are able to give a complete answer. $G$ is
 {\it{not}} relatively pseudo-Anosov if and only if 
$$((a_1,a_2), \{n_1,n_2\})
=(1,\{1,1\}),(1,\{1,2\}),(1,\{1,3\}),  {\text{ or }}  (2, \{1,1\})$$
(see Theorem 2.10).
 
In the case of $h \ge 3$ we get some bounds on the powers of the generators
sufficient for the group $G$ to be relatively pseudo-Anosov in the case 
each pair of the curves fill up the surface (see Theorem 3.6).
 We think that this should be answered in a more satisfactory way. 
This could be the subject of a 
forthcoming paper.

\proclaim{Question 0.2} Under what conditions $\G= \< D_{a_1},
\cdots, D_{a_h}\>$ is relatively pseudo-Anosov?
\endproclaim

In $\S 1$ we go over basic fact about Dehn twists and geometric intersection 
pairing and different kinds of ping-pong arguments we are going to use.

In $\S 2$ we determine exactly when a group generated by powers of 2 Dehn twists is free or relatively pseudo-Anosov.

In $\S 3$ we generalize the arguments for the case of $h \ge 3$ powers of
Dehn twists.

$\S 4$ is devoted to the simplest case of a group generated by powers of 3 Dehn twists,
i.e., when the curves are the $(1,0)$, $(0,1)$, and the $(1,1)$ curve on
 the torus. This group was studied previously in a more general setting (see
[BM], [Sch]).
 We show how to apply our methods in this case, and moreover we determine 
which ones do not contain surprising parabolic elements.

{\remark{Acknowledgment} I thank Daryl Cooper, Darren Long and
Marty Scharlemann for all the insightful conversations.

\head \S 1. Basics \endhead

For two simple closed curves $a,b$ let $(a,b$) denote their geometric intersection number.
For a set of simple closed curves $A=\{\a_1,\cdots, \a_n\}$ and a
simple closed curve (or measured lamination) $x$  put
$$||x||_A= \sum_{i=1}^n (x,\a_i).$$

Let $D_a$ be the (right-handed) Dehn twist in curve $a$.
The following lemma is proved in [FLP].

\proclaim{Lemma 1.1} For simple closed curves $a,b,c$, and $n \ge 0$,
$$ |(D^{\pm n}_a (b),c)-n(a,b)(a,c)| \le (b,c). \ \spadesuit$$  
\endproclaim

\proclaim{Lemma 1.2}(The ping-pong lemma; PPL)
Let $G$ be a group generated by $f_1,\cdots, f_n$, $n\ge 2$. Suppose 
 $G$ acts on a set $Z$. Assume that 
 there are $n$  non-empty mutually disjoint subsets $X_1,\cdots,X_n$ of $Z$
such that $f_i^{\pm n}(\cup_{j \ne i} X_j) \subset X_i$, for all $1 \le i \le n$.
Then $G \cong \Bbb F_n$.
\endproclaim

{\demo{Proof}} First notice that reduced  words of the form
$w=f_1^*f_i^* \cdots f_j^* f_1^*$ (*'s are non-zero integers)
are not the identity because $w(X_2) \cap X_2= \emptyset$. But any non-trivial word
in $f_1^*,\cdots,f_n^*$ is conjugate to a $w$ of the above form. $\spadesuit$.

\proclaim{Lemma 1.3}(ping-pong in the world trade center; PPWTC)
Let $G$ be a group generated by $f_1,\cdots, f_n$. Suppose 
 $G$ acts on a set $Z$, and there is a function  defined 
$||.||:Z \to \Bbb R_{\ge 0}$, with the following properties: 
 There are $n$  non-empty mutually disjoint subsets $X_1,\cdots,X_n$ of $Z$
such that $f_i^{\pm n}(Z-X_i) \subset X_i$ 
and for any $x \in Z-X_i$, we have $||f^{\pm n}(x)||>||x||$  for all
$n>0$.
Then $G \cong \Bbb F_n$. Moreover, the action of every  $g \in G$
which is not conjugate to some power of some $f_i$   on 
$Z$ has no 
periodic points.
\endproclaim

{\demo{Proof}}
Any reduced  word in $f_1^*,\cdots,f_n^*$ (*'s denote non-zero integers)
is conjugate to a reduced word $w=f_1^*\cdots f_1^*$.
To show that $w \ne id$ notice that if $x_1 \in Z-X_1$, then
$w(x_1) \in X_1$, therefore $w(x_1) \ne x_1$.
To prove the last assertion, notice that it's enough to show 
the claim with ``fixed points'' replaced by ``periodic points''.
Any element of $G$ which is not conjugate to a power of some $f_i$
is conjugate to some word of the form $w=f_j^* \cdots f_i^*$
with $i \ne j$. Now suppose $w(x)=x$. First assume 
$x \in Z-X_i$. Then by assumption $||w(x)||>||x||$
which is impossible. If on the other hand, $x \in X_i$ and 
$w(x)=x$, then $w\i(x)=f_i^* \cdots f_j^*(x)=x$. But again by assumption
$||w\i(x)||>||x||$, which is a contradiction.
  $\spadesuit$

Sometimes it is hard or  impossible to come up with subsets $X_1,\cdots,X_n$
Satisfying the assumptions of PPWTC. Here is a weaker version we are going to
use.

\proclaim{Lemma 1.4}(Weak ping-pong; WPP)
Let $G$ be a group generated by $f_1,\cdots, f_n$. Suppose 
 $G$ acts on a set $Z_1 \supset Z$, and there is a function  defined 
$||.||:Z \to \Bbb R_{\ge 0}$, with the following properties: 
 There are $n$  non-empty mutually disjoint subsets $X_1,\cdots,X_n$ of $Z$
such that $f_i^{\pm n}(Z-X_i) \subset X_i$ 
and for any $x \in Z-X_i$, we have $||f^{\pm n}(x)||>||x||$  for all
$n>0$. Also, assume that there exists an $n_0$ such that for every irreducible word  
$w=f_{k_1}^*\cdots f_{k_{n_0}}^*$,  we have $w(Z_1-Z) \subset X_{k_1}$.
Then $G \cong \Bbb F_n$. Moreover, the action of every  $g \in G$
which is not conjugate to some power of some $f_i$   on 
$Z$ has no 
periodic points. Also, if there is some $f \in G$ and $x \in Z_1-Z$
such that $f^n(x)=x$, then $f=w^k$ for some word $w$ of length $\le n_0$ and
$k \ge 1$. 
\endproclaim

\demo{Proof}
Easily obtained by altering the proof of 1.3. $\spadesuit$
\enddemo

\remark{Remark 1.5} In application of WPP, one must notice that,
 under the assumptions of the lemma, one may replace the assumption
$w(Z_1-Z) \subset X_{k_1}$ with $w_j(Z_1-Z) \subset X_{k_j}$ 
for some subword of $w$ of the form
 $w_j=f_{k_j}^*\cdots f_{k_{n_0}}^*$.

\proclaim{Lemma 1.6}(Cauchy-Schwartz inequality) If simple closed curves 
 $A=\{\a_1,\cdots, \a_n\}$ fill up $S$, then for any two 
simple closed curves $c_1,c_2$,
\roster
\item If $S$ is a closed surface, then
$$ (c_1,c_2) \le ||c_1||_A ||c_2||_A,$$
\item If $S$ has punctures
$$ (c_1,c_2) \le 2  ||c_1||_A ||c_2||_A.$$
\endroster

\endproclaim

{\demo{Proof}}
$S-\cup A$ is a union of discs or once punctured discs. 
One can put $c_1$ and $c_2$ in a transversal position 
such that each  segment of $c_1 - \cup A $ and $c_2 - \cup A $
intersect at most once in non-punctured discs, and intersect at most twice in
the punctured disks, while no three curves in $A \cup \{c_1,c_2\}$
pass through the same point. Therefore for each pair of segments in 
$ c_1 - (\cup A)$ and $ c_2 - (\cup A)$ one gets at most 
1 (resp. 2) intersection point(s) for closed (resp. punctured) $S$, 
hence the inequality. $\spadesuit$

\head \S 2. Groups generated by powers of two Dehn twists \endhead

Let $a,b$ be simple closed curves on $S$  with $(a,b)=m>0$. Put $A=\{a,b\}$.
Consider the sets of isotopy classes of simple closed curves
$$X=\{x \ | \ ||x||_A >0 \ \}.$$
Then for $\l \in (0,\infty)$ set
$$N_{a,\l}=\{x \in X \ | \  \ (x,a)<\l (x,b)\},$$
$$N_{b, \l\i}=\{x \in X\ | \  \ \l (x,b)<(x,a)\}.$$
Notice that $a \in N_{a,\l}$ and $b \in N_{b, \l\i}$,
 and $N_{a,\l} \cap N_{b,\l\i}=\emptyset$.

\proclaim{Lemma 2.1} 
\roster
\item $D_a^{\pm n}(N_{b,\l\i}) \subset N_{a,\l}$ if $mn \ge 2\l\i$.
\item If $mn \ge 2\l\i$ and $x \in  N_{b,\l\i}$, then
$|| D_a^{\pm n}(x)||_A > ||x||_A.$
\item $D_b^{\pm n}(N_{a,\l}) \subset N_{b,\l\i}$ if $mn \ge 2\l$. 
\item If $mn \ge 2\l$ and $x \in  N_{a,\l}$, then
$|| D_b^{\pm n}(x)||_A > ||x||_A.$
\endroster
\endproclaim

\demo{Proof}
Suppose $x \in N_b$ and $n>0$. Then

$$\align (D^{\pm n}_a(x), b) &\ge mn(x,a)-(x,b) \\
			&> (mn-\l\i)(x,a) \\
			&=(mn-\l\i)(D^{\pm n}_a(x), a).
\endalign$$
This shows that $D^{\pm n}_a(N_b) \subset N_a$ if $mn \ge 2\l\i$, which proves (1).
 By symmetry we immediately get (3). Now notice that for $x \in N_b$,
$$\align ||D^{\pm n }_a(x)||_A &=(D^{\pm n }_a(x),a)+(D^{\pm n }_a(x),b) \\
			&\ge (x,a)+ mn(x,a)-(x,b) \\
			&> (1+mn-\l\i)(x,a) \\
			&= (1+mn-\l\i)(1+\l)\i((x,a)+\l(x,a)) \\
			&>(1+mn-\l\i)(1+\l)\i(\l(x,b)+\l(x,a)) \\
			&=\l(1+mn-\l\i)(1+\l)\i||x||_A. 
\endalign$$
But $\l(1+mn-\l\i)(1+\l)\i \ge 1$ iff $mn \ge 2\l\i$, which proves (2), and 
by symmetry (4). $\spadesuit$
\enddemo

\proclaim{Theorem 2.2} For Two simple closed curves $a,b$ on a surface $S$,
\roster
\item If $(a,b)\ge 2$, the group $\langle D_a, D_b
 \rangle \cong \Bbb F_2$ 
(the free group on two generators).
\item If $(a,b)=1$, the group $\langle D_a^{n_1}, D_b^{n_2} \rangle \cong \Bbb F_2$ for $n_1,n_2 \ge 2$. 
\item
 If $(a,b)=1$ then $ \langle D_a, D_b^n \rangle \cong \Bbb F_2$ if $n\ge 4$.
\endroster
\endproclaim
 
\demo{Proof}
First consider the case $\l=1$. The group generated by $D_a$ and $D_b$ acts on $X$, and if $nm \ge 2$ one can apply PPL. But if
$m \ge 2$ this condition is automatically satisfied, and if $m=1$ one can apply it if both exponents are $\ge 2$. This proves (1),(2). 

To prove (3), consider the case $\l=2$. Lemma 2.1 shows that
 $D_a^{\pm n}(N_{b,1/2}) \subset N_{a,2}$ if $n \ge 1$, and 
 $D_b^{\pm n}(N_{a,2}) \subset N_{b,1/2}$ if $n \ge 4$. This finishes the 
proof by PPL. $\spadesuit.$
\enddemo

\proclaim{Corollary 2.3} For two simple closed curves $a,b$, the group 
$\langle D_a, D_b \rangle \cong \Bbb F_2 $ iff $(a,b) \ge 2$.
\endproclaim

{\demo{Proof}} If $(a,b)=0$,the twists $D_a, D_b$ commute. If $(a,b)=1$ then
the twists $D_a, D_b$ satisfy the braid relation
$D_aD_bD_a=D_bD_aD_b$. If $(a,b)\ge 2$, use Theorem 2.2. $\spadesuit$

\proclaim{Theorem 2.4}
Let $A=\{a,b\}$ be a set of two
 simple closed curves on a surface $S$ and $n_1,n_2>0$ be integers.
Put $G=\langle D_a^{n_1}, D_b^{n_2}\rangle$. The following conditions 
are equivalent:
\roster
\item
$G \cong \Bbb F_2$.
\item
Either $(a,b)\ge 2$, or 
$(a,b)=1$  and 
$$\{n_1,n_2\} \notin \{ \{1,1\}, \{1,2\}, \{1,3\} \}.$$
\endroster
\endproclaim

\demo{Proof}
We already saw that (2) implies (1). To prove (1) implies (2),
we must show that for $(a,b)=1$, the groups  $ \langle D_a, D_b^n \rangle$
 are not free for $n=1,2,3$. 
Let's denote $D_a$ by $a$ and $D_b$ by $b$ for
brevity. The case $n=1$ is non-free because $aba=bab$.
Now consider the case $n=2$. Observe that 
$$(ab^2)^2=ab^2 a b^2= ab(bab)b=ab(aba)b=(ab)^3,$$
so
$(ab^2)^4 a =(ab)^6a =a(bab)(aba)(bab)(aba)=a(aba)(bab)(aba)(bab)=a (ab)^6
=a (ab^2)^4.$
Which gives the  relation $(ab^2)^4a= a(ab^2)^4$   in $\langle a,b^2 \rangle$. 
In the case $n=3$, notice that
$(ab^3)^3=ab^3ab^3ab^3=ab^2(bab)b(bab)b^2=ab^2abababab^2= 
ab(bab)(aba)(bab)b=ab(aba)(bab)(aba)b=(ab)^6.$
Therefore we get the relation $(ab^3)^3a =a(ab^3)^3$.
$\spadesuit$

In the rest of this section we try to make sense out of the question ``which elements in $\langle D_a, D_b \rangle
$ are pseudo-Anosov?''

\proclaim{Theorem 2.5} Let $a,b$ be simple closed curves. Put $A=\{a,b\}$.
\roster
\item If $(a,b) \ge 3$, then the group $\langle D_a, D_b \rangle$
is pure. Moreover, each element is either conjugate to a power of 
$D_a$ or $D_b$, or restricts to a pseudo-Anosov element of $S_A$, and 
 is the identity on $S-S_A$.
\item If $(a,b)=2$ then the same conclusion of (1) holds for 
$\langle D_a, D_b^n \rangle$ and $\langle D_a^n, D_b \rangle$ 
for $n \ge 2$.
\item If $(a,b)=1$ then the same conclusion of (1) holds for 
$\langle D_a^{n_1}, D_b^{n_2} \rangle$ and
 $\langle D_a^{n_2}, D_b^{n_1} \rangle$ 
for $n_1 \ge 2$ and $n_2 \ge 3$.
\item If $(a,b)=1$ then the same conclusion of (1) holds for 
$\langle D_a^{n_1}, D_b^{n_2} \rangle$ and
 $\langle D_a^{n_2}, D_b^{n_1} \rangle$ 
for $n_1 \ge 1$ and $n_2 \ge 5$.
\endroster
\endproclaim

An element  $f \in M(S)$ is called {\it{pure}} (cf. [I])
if for any simple closed curve
$c$, $f^n(c)=c$ implies $f(c)=c$. In other words, by Thurston classification,
there is a finite
(possibly empty) set $C$ of disjoint  simple closed curves such that
$f$ fixes all curves in $C$, all components of $S- \cup C$, and is either
identity or pseudo-Anosov on each such component. A subgroup of $M(S)$ is 
called pure if all elements of it are pure. Ivanov showed that $M(S)$ 
is virtually pure, i.e., it has pure subgroups of finite index. Namely,
$ker(M(S) \to H_1(S, m \Bbb Z))$ for $m \ge 3$.

\demo{Proof}
Let $\l \in (1,3/2)$ be an irrational number. Then $X= N_{a,\l} \cup 
N_{b, \l\i}$. Now if $mn \ge 2\l\i$ Parts (1) and (2) of Lemma 2.1 
are satisfied, and if $mn \ge 2 \l$ parts (3),(4), and we can use PPWTC.
 But if $m \ge 3$, $n$ can be as low as 1 in both cases.
If $m=2$, $n_2$ has to be at least 2.
 If $m=1$, we must have $n_1 \ge 2$ and $n_2 \ge 3$.
 PPWTC shows that
no element in the corresponding groups can have a periodic curve which 
intersects $a$ or $b$. This means that all such elements are pseudo-Anosov
on $S_A$. On the other hand, all elements  leave all the points fixed 
outside a regular neighborhood of $a \cup b$.
To prove (4), take $\l \in (2,5/2)$ an irrational number. The first condition 
is  $n_1 \ge 1$ and the second one is $n_2 \ge 5$.
$\spadesuit$

{\proclaim{Corollary 2.6}} In Theorem 2.5, the groups in each case are 
relatively pseudo-Anosov.

\remark{2.7 Remark} Most of these pure groups are not 
contained in the ones discovered by
Ivanov, i.e., of the form $ker(M(S) \to H_1(S, t \Bbb Z))$, $t \ge 3$, if at least one of the simple closed curves $a,b$ is non-separating.
\endremark

This in particular proves

\proclaim{Corollary 2.8}(Thurston; see [FLP] etc.)
If $a,b$ fill up the closed surface $S=S(g)$, $g \ge 2$
 then $\langle D_a,D_b \rangle$
is free and all elements not conjugate to the powers of 
$D_a$ and $D_b$ are pseudo-Anosov.
\endproclaim
\demo{Proof} If $a,b$ fill up $S$ we must have $(a,b)\ge 3$. $\spadesuit$
\enddemo

\remark{2.9 Remarks} 

(i) The non-free cases $\<D_a,D_b\>$, $\<D_a, D_b^2\>$ and
$\<D_a, D_b^3\>$ where $(a,b)=1$ are not relatively  pseudo-Anosov because 
the maps $(D_aD_b)^6$, $(D_aD_b^2)^4$ and $(D_aD_b^3)^3$ commute with $D_a$
(see proof of 2.4), so they all fix $a$, and hence are reducible on $S_A$.

(ii)
If $(a,b)=2$ then $G=\<D_a,D_b\>$ is not relatively pseudo-Anosov. To 
prove this we consider two cases. The first case is when $S_A$ is a twice 
or once punctured torus. In this case  both $a,b$ 
can be embedded in a punctured torus
subsurface of $S$, say $a=(1,0)$ and $b=(1,2)$. It is easily seen by
geometric inspection that $D_bD_a$ fixes the curve $(1,1)$ and so is not
pseudo-Anosov. Case 2 is when  $S_A$
is a 4-punctured sphere. Call the punctures 1,2,3,4.
  Let $a$ be the curve bounding the punctures 
1,2, and $b$ be the curve bounding the punctures 2,3.
Then $D_a\i D_b$ fixes the curve bounding the punctures 1,3.

(iii) If $(a,b)=1$ then $G=\< D_a^2, D_b^2 \>$ is not relatively pseudo-Anosov.
For $D_b^2D_a^2$ is reducible.

(iv) If $(a,b)=1$ then $G=\<D_a,D_b^4\>$ is relatively pseudo-Anosov. Notice that this is the only case which the sets $N_{a,\l}$ and 
$N_{b,\l\i}$ and the PPWTC fail to show this fact, and we utilize WPP
in this case. One
 can assume  without loss of generality that $a=(1,0)$ and $b=(0,1)$ on
 a punctured torus. The only simple closed curves which are not 
in $X=N_{a,1} \cup N_{b,1}$ are $(1,1)$ and $(1,-1)$. But both of these
curves are mapped in $N_{a,1}$, $N_{b,1}$
 under $D_a^n$ ($n \ge 1$) or $D_b^k$, respectively ($k \ge 4$).
This proves that they can not be periodic points for any element of $G$, using
WPP with $n_0=1$.

We have shown

\proclaim{Theorem 2.10}
Let $A=\{a,b\}$ be a set of two
 simple closed curves on a surface $S$ and $n_1,n_2>0$ be integers.
Put $G=\langle D_a^{n_1}, D_b^{n_2}\rangle$. The following conditions 
are equivalent:
\roster
\item
$G$ is relatively pseudo-Anosov.
\item
Either $(a,b)\ge 3$, or $(a,b)=2$ and $(n_1,n_2)\ne (1,1)$, or
$(a,b)=1$  and 
$$\{n_1,n_2\} \notin \{ \{1,1\}, \{1,2\}, \{1,3\}, \{2,2\}  \}. \ \spadesuit$$
\endroster
\endproclaim

\head $\S 3.$ Groups generated by $h \ge 3$ powers of twists \endhead

In this section the phrase ``$i\ne j \ne k$'' means that $i,j,k$ are distinct.
Let $a_1,\cdots a_h$ be $h \ge 3$ simple closed curves on a surface $S$ such
 that $(a_i,a_j)>0$ for $i\ne j$.

Let  $\l_{ijk}>1$ and $\mu_{ij}>0$ (for $i\ne j \ne k$)
be real numbers such that $\mu_{ji}=\mu_{ij}\i$. Put $\l=(\l_{ijk})_{i\ne j \ne k}$ and $\mu=(\mu_{ij})_{i \ne j}$. Define the set of simple closed curves

$$N_{a_i}=N_{a_i,\l, \mu}=
\{x \ | \ (x,a_i)<\mu_{ij}(x,a_j), \ {(x,a_k) \over (x,a_j)}<\l_{ijk} {
(a_i,a_k) \over (a_i,a_j)},\
\forall j \ne k \ne i \},$$

for $i=1,\cdots, h$. 

\proclaim{Lemma 3.1} Let $a_1,\cdots,\a_h$ be a set of simple closed curves such that
$(a_i,a_j)\ne 0$ for $i \ne j$.
\roster 
\item For any choice of $\l_{ijk} \ge 1$, the sets $N_{a_i}$,
$i=1,\cdots,h$ are mutually disjoint.
\item
 For $1 \le i\ne j \le  h$, we have $D^{\pm n}_{a_i}(N_{a_j})
\subset  N_{a_i}$ for
$$\align n \ge \max \{& { 2 \over \mu_{ij}(a_i,a_j)}, \\
&  {1 \over \mu_{ik}(a_i,a_k)} + \l_{jik} {(a_j,a_k) \over (a_i,a_j)(a_i,a_k)},
\\
& { \l_{jil} \over \l_{ikl}-1}{(a_j,a_l) \over(a_i,a_l)
 (a_j,a_i)}+{\l_{ikl} \l_{jik} \over \l_{ikl}-1}
{(a_j,a_k)
\over (a_j,a_i)(a_i, a_k)}, \\
& { 1 \over (\l_{ikj}-1)\mu_{ij}(a_i,a_j)}+{\l_{ikj} \l_{jik} \over 
\l_{ikj}-1}
{(a_j,a_k)
\over (a_j,a_i)(a_i, a_k)}, \\
&  {\l_{ijl} \over (\l_{ijl}-1)\mu_{ij}(a_i,a_j)}
+{\l_{jil} \over \l_{ijl}-1}{(a_j,a_l) \over
(a_j,a_i)(a_i,a_l)}
\}_{k\ne l\ne i}.
\endalign$$
\endroster
\endproclaim

\demo{Proof} (1) is clear. To prove (2),
consider $x \in N_{a_j}$. We have
$$(D^{\pm n}_{a_i}(x),a_j) \ge n(a_i,a_j)(x,a_i)-(x,a_j)
> \mu_{ji}(x,a_i)$$
for $n \ge{ 2 \over \mu_{ij}(a_i,a_j)}$.
Let $k\ne i,j$. Then
$$(D^{\pm n}_{a_i}(x),a_k) \ge n(a_i,a_k)(x,a_i)-(x,a_k)
>\mu_{ki} (x,a_i)$$
if 
$$ n \ge {1 \over \mu_{ik}(a_i,a_k)} + \l_{jik} {(a_j,a_k) \over (a_i,a_j)(a_i,a_k)}.$$
Let $k,l\ne i$. Then
$$(a_l, D^{\pm n}_{a_i}(x))/(a_k, D^{\pm n}_{a_i}(x))<\l_{ikl}
(a_l,a_i)/(a_k, a_i)$$
iff
$$ {   (a_i, a_k)(D_{a_i}^{\pm n}(x), a_l) }
< \l_{ikl} {(a_i,a_l)  (D_{a_i}^{\pm n}(x), a_k)}.$$
This will hold if 
$${ (a_i, a_k)(n(a_i,a_l)(x,a_i)+(x,a_l)) } <
 \l_{ikl} {(a_i,a_l)(  n(a_i,a_k)(x,a_i)-(x,a_k)) }. \tag{*}$$
The inequality (*) is equivalent to 
$$n(a_i,a_l)(\l_{ikl}-1) > {(x,a_l) \over (x,a_i)}+\l_{ikl} {(x,a_k)(a_i,a_l)
\over (x,a_i)(a_i, a_k)}.\tag{**}$$
(one has
$(x,a_i)>0$ since $x \in N_{a_j}$.)
Therefore for $l\ne j$ and $k\ne j$ it's enough to have
$$n(a_i,a_l)(\l_{ikl}-1) \ge \l_{jil} {(a_j,a_l) \over (a_j,a_i)}+\l_{ikl} 
\l_{jik} 
{(a_j,a_k)(a_i,a_l)
\over (a_j,a_i)(a_i, a_k)};$$
i.e.,
$$n\ge{ \l_{jil} \over \l_{ikl}-1}{(a_j,a_l) \over(a_i,a_l)
 (a_j,a_i)}+{\l_{ikl} \l_{jik} \over \l_{ikl}-1}
{(a_j,a_k)
\over (a_j,a_i)(a_i, a_k)}.$$
If $l=j$ (and so $ k\ne j$) then one can replace (**) with 
$$n(a_i,a_l)(\l_{ikl}-1) \ge \mu_{ji}+\l_{ikl} {(x,a_k)(a_i,a_l)
\over (x,a_i)(a_i, a_k)}$$
which gives 
$$n \ge{ 1 \over (\l_{ikj}-1)\mu_{ij}(a_i,a_j)}+{\l_{ikj} \l_{jik} \over 
\l_{ikj}-1}
{(a_j,a_k)
\over (a_j,a_i)(a_i, a_k)}.$$
If $k=j$ (and so $l\ne j$) one similarly needs
$$n \ge {\l_{ijl} \over (\l_{ijl}-1)\mu_{ij}(a_i,a_j)}
+{\l_{jil} \over \l_{ijl}-1}{(a_j,a_l) \over
(a_j,a_i)(a_i,a_l)}.
\spadesuit$$
\enddemo

This lemma conveys the idea that if the set $\{(a_i,a_j)\}_{i\ne j}$ is not
``too spread around'' then the group $\G=\<D_{a_1},\cdots, D_{a_h}\>$
is free on $n$ generators, as follows:

\proclaim{Theorem 3.2}
Let $a_1,\cdots,a_h$ be simple closed curves on a surface $S$ such that
$M \le m^2/6$ where $M=\max\{(a_i,a_j)\}_{i\ne j}$ and $m=\min 
\{(a_i,a_j)\}_{i\ne j}$. Then  $\G=\<D_{a_1},\cdots, D_{a_h}\> \cong \Bbb F_h$.
More generally, suppose that for all $i \ne j \ne k$ we have
$${(a_i,a_k) \over (a_i, a_j)(a_j,a_k)} \le {1 \over 6}.$$
Then the same conclusion holds.
\endproclaim

\demo{Proof} Put $\mu_{ij}=1$ and $\l_{ijk}=2$ in Lemma 3.1.
By assumption, for all $i\ne j \ne k$,
$${(a_i,a_k) \over (a_i, a_j)(a_j,a_k)} \le {1 \over 6}.$$
This implies $(a_i,a_j) \ge 6$ for all $i\ne j$, since otherwise it is impossible for
both of 
$${(a_i,a_k) \over (a_i, a_j)(a_j,a_k)} \  \text{and} \
{(a_j,a_k) \over (a_i, a_j)(a_i,a_k)}$$
to be $\le 1/6$.
Therefore, it is easily seen that $n_i=1$ satisfies the requirements of 
Lemma 3.1. $\spadesuit$
\enddemo

On the other hand, if one allows the set  $\{(a_i,a_j)\}_{i\ne j}$ to be 
``more spread around'' then one gets weaker results:

\proclaim{Theorem 3.3}
Let $a_1,\cdots,a_h$ be simple closed curves on a surface $S$ such that
$$M_0=\max \{ {(a_i,a_k) \over (a_i, a_j)(a_j,a_k)} \}_{i,j,k}\ge {1 \over 6}$$
 Then  $\G=\<D^n_{a_1},\cdots,
 D^n_{a_h}\> \cong \Bbb F_h$ for $n \ge 6M_0$.
\endproclaim

One can easily construct infinitely many examples of such groups 
which are non-free, as follows. For a set of simple closed 
curves $A=\{a_1,\cdots,a_h\}$ define the twist set of $A$ as
 $$T(A)=\{ g(a) \ | \ g \in \< D_{\a} \>_{\a \in A}, a \in A \}.$$
Now let $A=\{a,b\}$ with $(a,b)=n\ge 2$ 
and $c_1 \in T(A)$, and pick $c$ such that $(c,c_1)=1$.
(To be able to do this we need $c_1$ to be non-separating which is 
possible if at least one of $a,b$, say $a$ is so.) Then We claim that 
$\<D_a,D_b,D_c\> \ncong \Bbb F_3$. The reason is that
since $c_1=g(a)$ where $g \in \<D_a, D_b \>$, $D_{c_1}=gD_ag\i$. Also
$D_c,D_{c_1}$ satisfy the braid relation, which gives a relation in 
$\<D_a,D_b,D_c\>$. This, together with Theorem 3.2 above implies that
the set of intersection number is ``spread around''.

Put $A=\{a_1,\cdots,a_h\}$, and 
$$X=\{x \ s.c.c. \ | \ ||x||_A >0 \ \}.$$

A deficiency of the sets $N_{a_i, \l, \mu}$ is that they don't cover $X$
in general. On the other hand, if we set $\l_{ijk}=\infty$, for irrational
 $\mu_{ij}$ we get the sets 
$$N_{a_i,\mu}=N_{a_i,\infty, \mu}=\{ x \ | \ (x, a_i)<\mu_{ij} (x,a_j) \}$$
which gives a disjoint cover of $X$, and we can hope to use PPWTC.

Now we try to give some conditions on the $a_i$ and $\l$ in order to have 
$N_{a_i, \infty, \mu}=N_{a_i,\l,\mu}$.

\proclaim{Lemma 3.4} Suppose for any $i\ne j$, $S_{\{a_i, a_j \}}=S$.
Then $N_{a_i, \infty, \mu}=N_{a_i,\l,\mu}$ for 
$$\l_{ijk}=2(a_i,a_j)(1 +{(a_k,a_j) \over (a_i,a_k)}) 
 (\mu_{ij}+1).$$
\endproclaim

\demo{Proof}
We only have to show  $N_{a_i, \infty, \mu} \subset N_{a_i,\l,\mu}$.
(The other inclusion is clear.) Suppose $x \in N_{a_i, \infty, \mu}$.
Then for $i \ne j \ne k$ we have (using Cauchy-Schwartz)
$${(x,a_k) \over (x,a_j)} \le {2((a_i,a_k)+(a_k,a_j)) ((x, a_i)+(x,a_j))
\over (x,a_j) }< 2((a_i,a_k)+(a_k,a_j))  (\mu_{ij}+1).$$
Therefore the lemma follows. $\spadesuit$  
\enddemo

\proclaim{Lemma 3.5} For $x \in N_{a_j, \l, \mu}$ and  $i \ne j$  we have
$||D^{\pm n}_{a_i}(x)||_A >||x||_A$ if 
$$n \ge { 2 \over ||a_i||_A}(\mu_{ji}+ \sum_{k \ne i,j} \l_{jik} {(a_j,a_k)
\over (a_j,a_i)} )$$
\endproclaim

\demo{Proof}
We have $\sum_k (D^{\pm n}_{a_i}(x),a_k) > \sum_k (x, a_k)$ if
$$ n \sum_{k \ne i} (a_i,a_k)(x,a_i)> 2 \sum_{k \ne i} (x,a_k).$$
i.e.,
$$ n \sum_k (a_i,a_k) > 2 {(x,a_j) \over (x,a_i)}+2 \sum_{k \ne i,j}  {(x,a_k) \over (x,a_i)}.$$ 
The proof follows. $\spadesuit$
\enddemo

Putting Lemmas 3.1, 3.4 and 3.5 together we get:

\proclaim{Theorem 3.6}
Let $\{a_1,\cdots,a_h\}$ be a
 set of simple closed curves such that $S_{ \{a_i, a_j \}}=S$ for $i \ne j$.
 Then $G=\< D^{n}_{a_1},\cdots, D^{n}_{a_h}\>$ is 
relatively pseudo-Anosov if $n \ge \max \{ 6M/m, 4M/m+5 \}$, where
$m=\min I \ge 2$ and $M=\max I$  in which $I= \{ (a_i,a_j) | i
\ne j \}$. \endproclaim

\demo{Proof}
Put $\mu_{ij}=1$ and take $\l_{ijk}$ as given by Lemma 3.4. A fairly straightforward computation shows that the number $n$ given satisfies the requirement
of 3.5 and the first two of 3.1 $\spadesuit$

\enddemo

It would be nice to improve this result. See Question 0.2.

\head $\S 4$ Example: A group  generated by 3 twists \endhead

To show that how the sets $N_{a_j, \l, \mu}$ work, let's look at the case $a_1,a_2,a_3$  are simple closed curves
such that $(a_i,a_j)=1$ for $i\ne j$ on a torus.
 The group $G=G(n_1,n_2,n_3)=\< D_{a_1}^{n_1}, D_{a_2}^{n_2},D_{a_3}^{n_3} \>$
 was studied by Backmuth, Mochizuki, and 
Scharlemann ([BM], [Sch]) in the case of real exponents $n_i$. The result of 
[Sch] implies that $G$ is free on 3 generators if 
$${1 \over n_1}+{1 \over n_2}+{1 \over n_3} \le 1.$$
Here we show how the results of $\S 3$ imply this.
Without loss of generality we may assume $a_1=(1,0)$, $a_2=(0,1)$, and $a_3=(1,1)$. Then, we are studying the groups
$$G=\<\left(  \matrix  1 & -1 \\ 0 & 1  \endmatrix \right)^{n_1},
\left(  \matrix  1 & 0 \\ 1 & 1  \endmatrix \right)^{n_2},
\left(  \matrix  2 & -1 \\ 1 & 0  \endmatrix \right)^{n_3} \>$$

Since $a_i$ cut the torus into two triangles, 
one can easily see that for any simple closed curve $x$,
$$(x,a_i)=(x,a_j)+(x,a_k)$$
for some $\{i,j,k\}=\{1,2,3\}$.
Suppose the numbers $\mu_{ij}$ satisfy the inequalities
$\mu_{ji}+\mu_{ki} \ge 1$. 
Now if for example  $(x,a_1)<\mu_{12}(x,a_2)$ and 
$(x,a_1)<\mu_{13}(x,a_3)$ then
either $(x,a_2)=(x,a_1)+(x,a_3)$ which implies
$$ 1 \le {(x,a_2) \over (x,a_3)} <1+\mu_{13}$$
or,  $(x,a_3)=(x,a_1)+(x,a_2)$ which implies 
$$ 1 \le {(x,a_3) \over (x,a_2)} <1+\mu_{12}.$$
(One cannot have $(x,a_1)=(x,a_2)+(x,a_3)$ since this implies
$\mu_{21}+\mu_{31}<1$.)
This means that, $N_{a_i,\l,\mu}=N_{a_i, \infty, \mu}$ for
$$\l_{ijk}=1+\mu_{ij}.$$
So one need not worry about making sure the conditions 
$$(x,a_k)/(x,a_j)<\l_{ijk}(a_i,a_k)/(a_i,a_j)$$
are satisfied. This together with Lemma 3.1 shows that $G(n_1,n_2,n_3) \cong \Bbb F_3$ if 
$$ n_i \ge \max \{ 2\mu_{ji}, 1+\mu_{ji}+\mu_{ki} \}_{j \ne k  \ne i}.$$

Now putting $\mu_{ij}=1$ implies that $G \cong \Bbb F_3$ if 
$n_i \ge 3$ for all $i$.

 On the other hand, if we assume
$\mu_{21}=\mu_{31}=1/2$ and $\mu_{32}=1$ we get that $G$ is free
if $n_1=2$ and $n_2,n_3 \ge 4$. 
 
Lastly, assume $\mu_{21}=2/3$, $\mu_{31}=1/3$ and $\mu_{32}=1/2$.
Then one gets that $G$ is free for $n_1=2$, $n_2=3$ and $n_3 \ge 6$.

This shows that $G \cong \Bbb F_3$ if $1/n_1+1/n_2+1/n_3 \le 1$.

To see for what values of $n_i$ the group $G$ is relatively Anosov we could 
perturb the numbers $\mu_{21},\mu_{31}, \mu_{32}$ to be irrational:

\proclaim{Theorem 4.1} Let $G$ be the  group $\<D_{a_1}^{n_1}, D_{a_2}^{n_2}, D_{a_3}^{n_3}
 \>$
where $a_i$ are simple closed curves on a torus or punctured torus  
with $(a_i,a_j)=1$ for $i \ne j$ and $n_i$ are positive integers.
\roster
\item $G \cong  \Bbb F_3$ if $1/{n_1}+ 1/{n_2}+1/n_3 \le 1$.
\item $G$ is relatively Anosov if $1/{n_1}+ 1/{n_2}+1/n_3 <1 $.
\endroster
\endproclaim

\demo{Proof}
(1) was already proved. To prove (2), notice that
the situation is completely symmetric 
so one can assume that $a=a_1=(1,0)$, $b=a_2=(0,1)$ and $c=a_3=(1,1)$.
PPWTC shows that for any choice of $\mu_{ij}$ such that all three of
$\mu_{21},\mu_{31},\mu_{32} $  are irrational and 
 $\mu_{ji}+\mu_{ki} \ge 1$ for $i \ne j \ne k$, $G$ is Anosov if
$$n_i \ge \max_{j} \{1+2\mu_{ji} \}=\max_{j \ne k} \{ 2 \mu_{ji}, 1+\mu_{ji}+
\mu_{ki}, 1+2\mu_{ji} \}.$$

Therefore, if $\mu_{21}=\mu_{31}=1-\e$ and $\mu_{32}=1+\e$
 where $\e$ is a small 
irrational, then we get the condition $n_1 \ge 3$ and $n_2,n_3 \ge 4$.

Similarly if
$\mu_{21}=2/3+\e$ and $\mu_{31}=1/3-\e$ and $\mu_{32}=1/2+\e$,
 then we get the condition $n_1 \ge 2$ and $n_2 \ge 4$ and $n_3 \ge 8$.

This means that $G$ is not relatively pseudo-Anosov only possibly
When 
$$\{n_1,n_2,n_3\}=\{1,*,*\},\{2,2,*\},\{2,3,*\},\{3,3,*\}, \{2,4,n\}, \ 4\le n \le 7.$$ 

One can cut this list a little short by using WPP, as follows: With 
$\mu_{ij}=1$ the only curves not covered by $N_i$ are $(1,-1)$ and $(2,1)$.
Then one can actually check that for 
$(n_1,n_2,n_3)=(3,3,n)$, $n \ge 4$,  WPP can be applied with $n_0=3$.

Also, in the case $(n_1,n_2,n_3)=(2,3,n)$, $n >7$, considering 
$\mu_{21}=2/3$, $\mu_{31}=1/3$, and $\mu_{32}=1/2$, one can see 
that only the curves $(2,3)$, $(2,-3)$, $(4,3)$ are not covered 
by $N_i$. Now WPP can be applied with
$n_0=3$.

Similarly for  $(n_1,n_2,n_3)=(2,4,n)$, $n >4$, one can use WPP.

So the list of possible exceptions cuts down to
$$\{1,*,,*\},\{2,2,*\},\{2,3,4\},\{2,3,5\},\{2,3,6\},\{2,4,4\},\{3,3,3\}.$$ 
This finishes the proof of the Theorem. $\spadesuit$

\enddemo

{\remark{4.2. Remark}}
The number 1 is the best possible in the above theorem, since 
$D_b^4 D_a^2 D_c^4$ fixes $(1,2)$ and  $D_b^3 D_a^2 D_c^6$ fixes $(2,3)$. 
Also, $D_b^2 D_a^2 D_c^n$ fixes $(1,1)$, so this element commutes with
$D_c$.

\

{\bf{References}}
\roster

\item"[BM]" S. Bachmuth and H. Mochizuchi, Triples of $2 \times 2$ matrices which generate free groups, Proc. Am. Math Soc. 59, (1976) 25-28.

\item"[FLP]" A. Fathi, F. Laundenbach and V. Poenaru, Travaux de
Thurston sur les surfaces,
 Ast\'erisque 66-67,
 Soc. Math\'ematique de France,
 1979.

\item"[I]" N. Ivanov, Subgroups of Teichm\"uller modular groups,
 Translations of Mathematical Monographs, 115. Am. Math.
 Soc., Providence, RI, 1992.

\item"[Sch]"  M. Scharlemann, Subgroups of $SL(2, \Bbb R)$ freely generated 
by three parabolic elements, Linear and
Multilinear Algebra 7 (1979), no. 3, 177-191.

\endroster

\end

\proclaim{Lemma} Let $a,b,x$ be simple closed curves, and $n$ a non-zero integer. If
$$(D^{n}_a(x),b)=n(a,b)(x,a)-(x,b),$$
then
$$(D^{-n}_a(x),b)=n(a,b)(x,a)+(x,b).$$
\endproclaim

\demo{Proof}
We have
$$|(D_a^{-2n}(D_a^n(x)), b)-2n(a,b)(x,a)| \le (D^n_a(x),b).$$
Therefore,
$$ (D_a^{-n}(x), b) \ge 2n(a,b)(x,a)-(D^n_a(x),b)=n(a,b)(x,a)+(x,b). $$
On the other hand,
$$|(D_a^{-n}(x), b)-n(a,b)(x,a)| \le (x,b),$$
so
$$ (D_a^{-n}(x), b) \le n(a,b)(x,a)+(x,b),$$
which proves the lemma. $\spadesuit$
\enddemo